\input amstex
\documentstyle{amsppt}
\magnification=\magstephalf
\pagewidth{6.25truein}
\pageheight{9.25truein}
\parskip=10pt
\voffset=0pt
\font\smallit=cmti10
\font\smalltt=cmtt10
\font\smallrm=cmr9
\nologo
\topmatter
\leftheadtext{\hskip 8pt Zhi-Wei Sun}
\rightheadtext{\hskip 8pt On covering numbers}

\def\Z{\Bbb Z}

\def\l{\left}
\def\r{\right}
\def\bg{\bigg}
\def\({\bg(}
\def\[{\bg[}
\def\){\bg)}
\def\]{\bg]}
\def\t{\text}
\def\f{\frac}
\def\mo{\roman{mod}}

\def\sp {\supseteq}
\def\sm{\setminus}

\def\eq{\equiv}
\def\cs{\cdots}
\def\ls{\leqslant}
\def\gs{\geqslant}
\def\al{\alpha}

\def\Proof{\noindent{\it Proof}}

\def\Proof{\noindent{\it Proof}}

\def\Remark{\smallskip\noindent{\it  Remark}}
\def\Ack{\medskip\noindent {\bf Acknowledgment}}
\endtopmatter

\document \vskip -40pt
\centerline{\hskip 8pt To appear in
\smalltt INTEGERS: \smallrm ELECTRONIC JOURNAL OF COMBINATORIAL NUMBER THEORY}
\vskip 40pt

\centerline{\bf ON COVERING NUMBERS} \vskip 20pt \centerline{{\bf
Zhi-Wei Sun}\footnote {Supported by the National Natural Science
Fund for Distinguished Young Scholars (No. 10425103) and a Key
Program of NSF (No. 10331020) in China.}}
\medskip
\centerline {\smallit Department of Mathematics, Nanjing University,
Nanjing 210093, P. R. China}
\centerline{{\tt zwsun\@nju.edu.cn}\quad\ \ \   http://pweb.nju.edu.cn/zwsun/}
\vskip 20pt
\centerline{\it Dedicated to Prof. R. L. Graham for his 70th birthday}
\vskip 30pt
\centerline{\bf Abstract} A positive integer $n$ is called a
covering number if there are some distinct divisors
$n_1,\ldots,n_k$ of $n$ greater than one and some integers
$a_1,\ldots,a_k$ such that $\Z$ is the union of the residue
classes $a_1(\mo\ n_1),\ldots,a_k(\mo\ n_k)$. A covering number is
said to be primitive if none of its proper divisors is a covering
number. In this paper we give some sufficient conditions for $n$
to be a (primitive) covering number; in particular, we show that
for any $r=2,3,\ldots$ there are infinitely many primitive
covering numbers having exactly $r$ distinct prime divisors. In
1980 P. Erd\H os asked whether there are infinitely many positive
integers $n$ such that among the subsets of $D_n=\{d\gs 2:\, d\mid
n\}$ only $D_n$ can be the set of all the moduli in a cover of
$\Z$ with distinct moduli; we answer this question affirmatively.
We also conjecture that any primitive covering number must have a
prime factorization $p_1^{\al_1}\cdots p_r^{\al_r}$ (with
$p_1,\ldots,p_r$ in a suitable order) which satisfies
$\prod_{0<t<s}(\al_t+1)\gs p_s-1$ for each $1\ls s\ls r$, with
strict inequality when $s=r$.
\noindent \baselineskip=15pt


\vskip 30pt

\subhead \nofrills{\bf 1. Introduction}\endsubhead

For $a\in\Z$ and $n\in\Z^+=\{1,2,3,\ldots\}$,
$a(\mo\ n)=\{a+nx:\, x\in\Z\}$ is called a residue class
with modulus $n$.
If every integer lies in at least one of the residue classes
$a_1(\mo\ n_1), \ldots, a_k(\mo\ n_k),$
then we call the finite system
$$A=\{a_i(\mo\ n_i)\}^k_{i=1}\tag 1.0$$
a {\it cover} of $\Z$ (or {\it covering system}),
and $n_1,\cs,n_k$ its {\it moduli}.
If (1.0) forms a cover of $\Z$ but none of its proper subsystems does,
then (1.0) is said to be a {\it minimal cover} of $\Z$.

 In the 1930s P. Erd\H os (cf. [E50])
invented the concept of a cover of $\Z$ and gave the following example
$$\{0(\mo\ 2),\ 0(\mo\ 3),\ 1(\mo\ 4),\ 5(\mo\ 6),\ 7(\mo\ 12)\}$$
whose moduli $2,\,3,\,4,\,6,\,12$ are distinct.
Covers of $\Z$ with distinct moduli are of particular interest and they have
some surprising applications (see, e.g., [F] and [S00]).
For problems and results concerning covers of $\Z$ and their generalizations
the reader may consult [E97], [FFKPY], [Gu], [PS], [S03], [S04] and [S05].

 Here is a famous open conjecture.
\proclaim{The Erd\H os--Selfridge Conjecture} If $(1.0)$ forms a cover of $\Z$ with
the moduli $n_1,\ldots,n_k$ distinct and greater than one, then
 $n_1,\ldots,n_k$ are not all odd.
\endproclaim

Following J. A. Haight [H] we introduce the following concept.

\noindent{\bf Definition 1.1}. A positive integer $n$
is called a {\it covering number} if
there is a cover of $\Z$ with all the moduli distinct, greater than one and dividing $n$.

Erd\H os' example shows that $2^2\cdot3=12$ is a covering number.
By density considerations, if $n$ is a covering number then $\sum_{1<d\,\mid\,n}1/d\gs1$;
it follows that none of $2,3,\ldots,11$ is a covering number.
Moreover, Example 3 of [S96] indicates that $2^{n-1}n$
is a covering number for every $n=3,5,7,\ldots$.

In the direction of the Erd\H os--Selfridge conjecture,
S. Guo and Z. W. Sun [GS] proved that any odd and squarefree covering number
should have at least 22 distinct prime divisors.

If $(1.0)$ is a cover of $\Z$ with $n_1\ls\cdots\ls n_{k-1}<n_k$,
then $\sum_{i=1}^{k-1}1/n_i\gs1$ by Theorem I\,(iv) of
Sun [S96]. So, a necessary condition for $n\in\Z^+$ to be a
covering number is that
$$\f{\sigma(n)}n=\sum_{d\mid n}\f1d\gs2+\f1n,\tag1.1$$
where $\sigma(n)$ is the sum of all positive divisors of $n$.
However, as shown by Haight [H], there does not exist a constant $c>0$ such that
$n\in\Z^+$ is a covering number whenever $\sigma(n)/n>c$.

Let (1.0) be a cover of $\Z$, and set $w(r)=|\{1\ls i\ls k:\, r\eq a_i\ (\mo\ n_i)\}|$
for $r=0,\ldots,N-1$,
where $N=[n_1,\ldots,n_k]$ is the least common multiple of $n_1,\ldots,n_k.$
By Theorem 5(ii) and Example 6 of [S01],
$$\sum\Sb 1\ls i\ls k\\\gcd(x+a_i,n_i)=1\endSb\f1{\varphi(n_i)}
=\sum\Sb 0\ls r<N\\\gcd(x+r,N)=1\endSb\f{w(r)}{\varphi(N)}
\gs\sum\Sb 0\ls r<N\\\gcd(x+r,N)=1\endSb\f{1}{\varphi(N)}=1\ \quad\t{for all}\ x\in\Z,$$
where $\varphi$ is Euler's totient function.
If $1<n_1<\cdots<n_k$ and $x\eq -a_i\ (\mo\ n_i)$ for all those $i\in I=\{1\ls j\ls k:\,
n_j\ \t{is a prime}\}$ (such an integer $x$ exists by the Chinese Remainder Theorem), then
$$\sum^k\Sb i=1\\i\not\in I\endSb\f1{\varphi(n_i)}
\gs\sum\Sb 1\ls i\ls k\\\gcd(x+a_i,n_i)=1\endSb\f1{\varphi(n_i)}\gs1.$$
Thus, if $n\in\Z^+$ is a covering number then we have
$$\sum\Sb d\mid n\\ d\ \t{is composite}\endSb\f1{\varphi(d)}\gs1.\tag1.2$$

Throughout this paper, for a predicate $P$ we let
$[\![P]\!]$ be $1$ or $0$ according as $P$ holds or not.
For a real number $x$, as usual we use $\lfloor x\rfloor$ and $\lceil x\rceil$
to denote the greatest integer not exceeding $x$ and the least integer greater than or equal to $x$,
respectively.

Our first theorem in this paper gives a sufficient condition for covering numbers.

\proclaim{Theorem 1.1} Let $p_1,\ldots,p_r$ be distinct primes,
and let $\al_1,\ldots,\al_r\in\Z^+$.
Suppose that
$$\prod_{0<t<s}\l(\al_t+1\r)\gs p_s-[\![r\not=s]\!]\qquad\t{for all}\ s=1,\ldots,r.\tag1.3$$
Then $p_1^{\al_1}\cs p_r^{\al_r}$ is a covering number.
\endproclaim

\Remark\ 1.1. As usual the empty product $\prod_{0<t<1}\l(\al_t+1\r)$ is regarded as $1$,
thus (1.3) implies that $p_1=2\ls r$.

The Erd\H os--Selfridge conjecture can be viewed as the converse of the following result.

\proclaim{Corollary 1.1} Let $p_1=2<p_2<\cdots<p_r\ (r>1)$ be distinct primes.
Then there are $\al_1,\ldots,\al_r\in\Z^+$ such that
$p_1^{\al_1}\cdots p_r^{\al_r}$ is a covering number.
\endproclaim

\Proof. For $t=1,\ldots,r-1$ we set
$$\al_t=\l\lceil\f{p_{t+1}-[\![t\not=r-1]\!]}{p_t-1}\r\rceil-1.$$
Then
$$\prod_{0<t<s}(\al_t+1)\gs\prod_{0<t<s}\f{p_{t+1}-[\![t+1\not=r]\!]}{p_t-[\![t\not=r]\!]}=
\f{p_s-[\![s\not=r]\!]}{p_1-[\![1\not=r]\!]}=p_s-[\![r\not=s]\!]$$
for all $s=1,\ldots,r$. Thus, by Theorem 1.1, $p_1^{\al_1}\cdots p_r^{\al_r}$
is a covering number. \qed

In contrast with Corollary 1.1, we have the following second theorem.

\proclaim{Theorem 1.2} Let $\al_1,\ldots,\al_r\in\Z^+$. Then
$p_1^{\al_1}\cdots p_r^{\al_r}$ is a covering number for some distinct primes
$p_1<\cdots<p_r$,
if and only if one of the following {\rm (i)--(iii)} holds.

{\rm (i)} $r=2\ls\al_1;$\qquad\quad
{\rm (ii)} $r=3$ and $\max\{\al_1,\al_2\}\gs 2;$\qquad\quad
{\rm (iii)} $r\gs 4$.
\endproclaim

\noindent{\bf Definition 1.2}. A covering number is called a {\it primitive covering number}
if none of its proper divisors is a covering number.

Our third theorem provides a sufficient condition for
primitive covering numbers.

\proclaim{Theorem 1.3} Let $p_1=2<p_2<\cdots<p_r\ (r>1)$ be distinct primes.
Suppose further that $p_t-1\mid p_{t+1}-1$ for all $0<t<r-1$, and $p_r\gs(p_{r-1}-2)(p_{r-1}-3)$.
Then
$$p_1^{\f{p_2-1}{p_1-1}-1}\cdots p_{r-2}^{\f{p_{r-1}-1}{p_{r-2}-1}-1}
p_{r-1}^{\lfloor\f{p_r-1}{p_{r-1}-1}\rfloor} p_r$$
is a primitive covering number.
\endproclaim

\Remark\ 1.2. By Theorem 1.3, the number $2\cdot3\cdot5\cdot7=210$
is a primitive covering number; moreover, Erd\H os constructed a
cover of $\Z$ whose moduli are all the 14 proper divisors of $210$
(cf. [Gu] or [GS]).

\proclaim{Corollary 1.2} For any $r=2,3,\ldots$ there are infinitely many primitive
covering numbers having exactly $r$ distinct prime divisors.
\endproclaim
\Proof. By Dirichlet's theorem (cf. [R, pp.\,237--244]),
for any $m\in\Z^+$ there are infinitely many
primes $p$ such that $m\mid p-1$. So, the desired result follows from Theorem 1.3. \qed

As an application of Theorem 1.3 and its proof, here we give our last theorem.

\proclaim{Theorem 1.4}
{\rm (i)} An integer $n>1$ with at most two distinct prime divisors
is a primitive covering number if and only if
$n=2^{p-1}p$ for some odd prime $p$.

{\rm (ii)}  A positive integer $n\eq0\ (\mo\ 3)$ with exactly three distinct prime divisors
is a primitive covering number
if and only if $n=2\cdot3^{(p-1)/2}p$ for some prime $p>3$.

{\rm (iii)} If $p>5$ is a prime, then
both $2^35^{\lfloor(p-1)/4\rfloor}p$ and $2\cdot 3\cdot 5^{\lfloor(p-1)/4\rfloor}p$
are primitive covering numbers. If $p>7$ is a prime, then
$2\cdot 3^27^{\lfloor(p-1)/6\rfloor}p$ is a primitive covering number,
and so is $2^57^{\lfloor(p-1)/6\rfloor}p$ provided that $p\not=13,19$.
\endproclaim

\Remark\ 1.3. Note that $2^57^2\cdot 13$ and $2^57^3\cdot 19$ are both covering numbers
by Theorem 1.1. But we don't know whether they are primitive covering numbers.

The following corollary provides an affirmative answer to a question of Erd\H os [E80].

\proclaim{Corollary 1.3} There are infinitely many positive integers $n$ such that
among the subsets of $D_n=\{d\gs2:\, d\mid n\}$ only $D_n$ can be the set of all the moduli in a cover
of $\Z$ with distinct moduli.
\endproclaim
\Proof. Let $p$ be one of the infinitely many odd primes. By Theorem 1.4(i), $2^{p-1}p$
is a primitive covering number.

Let (1.0) be any minimal cover of $\Z$
with $1<n_1<\cdots<n_k$ and $[n_1,\ldots,n_k]=2^{p-1}p$.
We want to show that $\{n_1,\ldots,n_k\}=\{d>1:\, d\mid 2^{p-1}p\}$.
By a conjecture of \v S. Zn\'am proved by R. J. Simpson [Si], we have
$$k\gs 1+f([n_1,\ldots,n_k])=1+(p-1)(2-1)+(p-1)=2p-1,$$
where the Mycielski function $f:\Z^+\to\Z$ is given by
$f(\prod_{t=1}^rp_t^{\al_t})=\sum_{t=1}^r\al_t(p_t-1)$
with $p_1,\ldots,p_r$ distinct primes and $\al_1,\ldots,\al_r$ nonnegative integers
(cf. [S90] and [Z]).
On the other hand,
$$k\ls |\{d>1:\, d\mid 2^{p-1}p\}|=|\{2^\al p^{\beta}:\,
\al=0,\ldots,p-1;\ \beta=0,1\}|-1=2p-1.$$
So $k=2p-1=|\{d>1:\, d\mid 2^{p-1}p\}|$ and we are done. \qed

In the next section we are going to prove Theorems 1.1 and 1.2. Section 3
is devoted to our proofs of Theorems 1.3 and 1.4.
To conclude this section we propose the following conjecture concerning the converse of
Theorem 1.1.

\proclaim{Conjecture 1.1} Any primitive covering number can be written in the form
$p_1^{\al_1}\cdots p_r^{\al_r}$
with $p_1,\ldots,p_r$ distinct primes and $\al_1,\ldots,\al_r\in\Z^+$,
such that $(1.3)$ is satisfied.
\endproclaim

\Remark\ 1.4. Actually the author made this conjecture on July 16, 1988.
Since (1.3) implies $p_1=2$,
Conjecture 1.1 is stronger than the Erd\H os--Selfridge conjecture.

\subhead \nofrills{\bf 2. Proofs of Theorems 1.1 and 1.2}\endsubhead

For $n\in\Z^+$ let $d(n)$ denote the number of distinct positive divisors of $n$.
If $n$ has the factorization $p_1^{\al_1}\cdots p_r^{\al_r}$ where
$p_1,\ldots,p_r$ are distinct primes and $\al_1,\ldots,\al_r\in\Z^+$,
then it is well known that $d(n)=\prod_{t=1}^r(\al_t+1)$.

\noindent{\it Proof of Theorem 1.1}. For each $s=1,\ldots,r$, since
$$d(p_1^{\al_1}\cdots p_{s-1}^{\al_{s-1}})=\prod_{0<t<s}(\al_t+1)\gs p_s-[\![r\not=s]\!]$$
there exist $p_s-[\![r\not=s]\!]$
distinct positive divisors $d_1^{(s)},\ldots,d^{(s)}_{p_s-[\![r\not=s]\!]}$
of $\prod_{0<t<s}p_t^{\al_t}$.
Let $\Cal A$ be the system consisting of $0(\mo\ d^{(r)}_{p_r}p_r^{\al_r})$
and the following $\sum_{s=1}^r\al_s(p_s-1)$ residue classes:
$$jp_1^{\al_1}\cdots p_{s-1}^{\al_{s-1}}p_s^{\al-1}(\mo\ d_j^{(s)}p_s^{\al})
\ \ \ (\al=1,\ldots,\al_s;\ j=1,\ldots,p_s-1;\ s=1,\ldots,r).$$
Then all the moduli of $\Cal A$ are distinct. Observe that
$$\align&\bigcup_{j=1}^{p_s-1}jp_1^{\al_1}\cdots p_{s-1}^{\al_{s-1}}p_s^{\al-1}
(\mo\ d_j^{(s)}p_s^{\al})
\\\sp&\bigcup_{j=1}^{p_s-1}jp_1^{\al_1}\cdots p_{s-1}^{\al_{s-1}}p_s^{\al-1}
(\mo\ p_1^{\al_1}\cdots p_{s-1}^{\al_{s-1}}p_s^{\al})
\\=&0(\mo\ p_1^{\al_1}\cdots p_{s-1}^{\al_{s-1}}p_s^{\al-1})\sm
0(\mo\ p_1^{\al_1}\cdots p_{s-1}^{\al_{s-1}}p_s^{\al})
\endalign$$
and
$$\align&\bigcup_{\al=1}^{\al_s}\l(0(\mo\ p_1^{\al_1}\cdots p_{s-1}^{\al_{s-1}}p_s^{\al-1})
\sm 0(\mo\ p_1^{\al_1}\cdots p_{s-1}^{\al_{s-1}}p_s^{\al})\r)
\\&\quad=0(\mo\ p_1^{\al_1}\cdots p_{s-1}^{\al_{s-1}})
\sm 0(\mo\ p_1^{\al_1}\cdots p_{s-1}^{\al_{s-1}}p_s^{\al_s}).
\endalign$$
If an integer $x$ is not in the residue class
$0(\mo\ d^{(r)}_{p_r}p_r^{\al_r})$, then $x\not\eq0\ (\mo\ p_1^{\al_1}\cdots p_r^{\al_r})$
and hence
$$x\in0(\mo\ 1)\sm 0(\mo\ p_1^{\al_1}\cdots p_r^{\al_r})
=\bigcup_{s=1}^r\(0\(\mo\ \prod_{0<t<s}p_t^{\al_t}\)
\sm 0\(\mo\ \prod_{t=1}^s p_t^{\al_t}\)\).$$
Therefore $\Cal A$ does form a cover of $\Z$. \qed

\Remark\ 2.1. In the proof of Theorem 1.1, we make use of some basic ideas in [Z] and [S90].

\proclaim{Lemma 2.1}
Let $p_1,\ldots,p_r$ be distinct primes and $\al_1,\ldots,\al_r\in\Z^+$.
Suppose that $p_1^{\al_1}\cdots p_r^{\al_r}$ is a covering number
but $\prod_{0<t<r}p_t^{\al_t}$ is not.
Then we must have $\prod_{0<t<r}(\al_t+1)\gs p_r$.
\endproclaim
\Proof. Let (1.0) be a minimal cover of $\Z$
with $1<n_1<\cdots<n_k$ and $[n_1,\ldots,n_k]\mid p_1^{\al_1}\cdots p_r^{\al_r}$.
Since $\prod_{0<t<r}p_t^{\al_t}$ is not a covering number, $p_r$ divides
$[n_1,\ldots,n_k]$. Let $\al\in\Z^+$ be the largest integer such that
$p_r^\al$ divides at least one of the moduli $n_1,\ldots,n_k$.
Then we have
$$|\{1\ls i\ls k:\, p_r^{\al}\mid n_i\}|\gs p_r$$
by [SS, Theorem 1] or [S96, Corollary 3]. Note that
$$|\{1\ls i\ls k:\, p_r^{\al}\mid n_i\}|
\ls|\{dp_r^{\al}:\, d\mid p_1^{\al_1}\cdots p_{r-1}^{\al_{r-1}}\}|
=d\(\prod_{0<t<r}p_t^{\al_t}\)=\prod_{0<t<r}(\al_t+1).$$
So the desired result follows. \qed

\medskip
\noindent{\it Proof of Theorem 1.2}. If (i) holds,
then $2^{\al_1}3^{\al_2}$ is a covering number
by Theorem 1.1 since $1\gs 2-1$ and $\al_1+1\gs 3$. If (ii) is valid,
then $2^{\al_1}3^{\al_2}5^{\al_3}$ is a covering number
by Theorem 1.1, since $\al_1+1\gs 3-1$ and $(\al_1+1)(\al_2+1)\gs (1+1)(2+1)>5$.
When (iii) happens (i.e., $r\gs4$), letting $p_1,\ldots,p_r$ be the first $r$ primes
in the ascending order, we then have $p_1=2,\, p_2=3,\ p_3=5,\, p_4=7$,
hence
$\prod_{s=1}^rp_s^{\al_s}$ is a covering number
by Theorem 1.1, because $\al_1+1\gs 3-1$, $(\al_1+1)(\al_2+1)\gs 5-1$, and
 $p_s<2^{s-1}\ls\prod_{0<t<s}(\al_t+1)$ for $s\gs4$
(by mathematical induction and Bertrand's postulate (cf. [R, pp.\,220--221])
proved by Chebyshev).

Now suppose that there are distinct primes $p_1<\cdots<p_r$ such that
$n=p_1^{\al_1}\cdots p_r^{\al_r}$ is a covering number.
Let $d>1$ be the smallest covering number dividing $n$. Then
$d$ is a primitive covering number. By Lemma 2.1, $d$ cannot be a prime power.
So $r\gs 2$. If $r=2$ and $\al_1=1$,
then $d=p_1p_2^{\beta}$ for some $\beta=1,\ldots,\al_2$,
thus by Lemma 2.1 we get the contradiction
$1+1\gs p_2>p_1\gs 2.$
If $r=3$ and $\al_1=\al_2=1$, then
$d=p_1p_2p_3^{\gamma}$ for some $\gamma=1,\ldots,\al_3$,
hence by Lemma 2.1 we have
$(1+1)(1+1)\gs p_3\gs5$ which is impossible.
Therefore one of (i)--(iii) holds.
\qed

\subhead \nofrills{\bf 3. Proofs of Theorems 1.3 and 1.4}\endsubhead

\medskip
\noindent{\it Proof of Theorem 1.3}. Set
$$\al_1=\f{p_2-1}{p_1-1}-1,\ \ldots,\ \al_{r-2}=\f{p_{r-1}-1}{p_{r-2}-1}-1\ \
\t{and}\ \ \al_{r-1}=\l\lfloor\f{p_r-1}{p_{r-1}-1}\r\rfloor.$$
Then
$$\prod_{0<t<s}(\al_t+1)=\prod_{0<t<s}\f{p_{t+1}-1}{p_t-1}=\f{p_s-1}{p_1-1}=p_s-1$$
for $s=1,\ldots,r-1$, and
$$\prod_{0<t<r}(\al_t+1)
=\prod_{0<t<r-1}(\al_t+1)\times\(\l\lfloor\f{p_r-1}{p_{r-1}-1}\r\rfloor+1\)
>(p_{r-1}-1)\f{p_r-1}{p_{r-1}-1}=p_r-1.$$
Thus $n=p_1^{\al_1}\cdots p_{r-1}^{\al_{r-1}}p_r$ is a covering number
in light of Theorem 1.1.

Let $d>1$ be the smallest covering number dividing $n$.
It remains to show that $d=n$.

Suppose that $p_s$ is the maximal prime
divisor of $d$. If $s\not=r$, then
$\prod_{0<t<s}(\al_t+1)=p_s-1<p_s$
which contradicts Lemma 2.1. Therefore, $d$
has the form $p_1^{\beta_1}\cdots p_{r-1}^{\beta_{r-1}}p_r$
where $\beta_t\in\{0,\ldots,\al_t\}$
for $t=1,\ldots,r-1$. By Lemma 2.1,
$$\prod_{0<t<r}(\beta_t+1)\gs p_r.$$
If $\beta_{r-1}<\al_{r-1}$, then
$$\prod_{0<t<r}(\beta_t+1)\ls\prod_{0<t<r-1}(\al_t+1)\times\al_{r-1}
=(p_{r-1}-1)\l\lfloor\f{p_r-1}{p_{r-1}-1}\r\rfloor\ls p_r-1<p_r.$$
So we must have $\beta_{r-1}=\al_{r-1}$.

Assume that $\beta_j<\al_j$ for some $1\ls j\ls r-2$. Then
$$\prod_{t=1}^{r-1}(\beta_t+1)
\ls \prod_{t=1}^{r-2}(\al_t+1)\times\f{\al_j}{\al_j+1}(\al_{r-1}+1)=m,$$
where
$$\align m=&(p_{r-1}-1)\l(1-\f{p_j-1}{p_{j+1}-1}\r)
\(\l\lfloor\f{p_r-1}{p_{r-1}-1}\r\rfloor+1\)
\\\ls&(p_{r-1}-1)\l(1-\f{p_j-1}{p_{j+1}-1}\r)
\(\f{p_r-1}{p_{r-1}-1}+1\)
\\=&(p_{r-1}-2+p_r)\l(1-\f{p_j-1}{p_{j+1}-1}\r).
\endalign$$
Since
$$p_r\gs(p_{r-1}-3)(p_{r-1}-1-1)\gs(p_{r-1}-3)\(\f{p_{j+1}-1}{p_j-1}-1\),$$
we have
$$(p_{r-1}-2)\(\f{p_{j+1}-1}{p_j-1}-1\)-p_r<\f{p_{j+1}-1}{p_j-1}$$
and hence
$$m\ls (p_{r-1}-2)\(1-\f{p_j-1}{p_{j+1}-1}\)+p_r-p_r\f{p_j-1}{p_{j+1}-1}<p_r+1.$$
We claim that
$$m=(p_j-1)\l(\f{p_{r-1}-1}{p_j-1}-\f{p_{r-1}-1}{p_{j+1}-1}\r)
\(\l\lfloor\f{p_r-1}{p_{r-1}-1}\r\rfloor+1\)\not=p_r.$$
In fact, $m$ is composite when $j>1$; if $j=1$ then
$$\f{p_{r-1}-1}{p_j-1}-\f{p_{r-1}-1}{p_{j+1}-1}
=p_{r-1}-1-\f{p_{r-1}-1}{p_2-1}\gs\f{p_{r-1}-1}2>1$$
unless $p_{r-1}=3$ in which case
$$m=\l\lfloor\f{p_r-1}{3-1}\r\rfloor+1=\f{p_r+1}2<p_r.$$
In view of the above,
$$p_r\ls \prod_{t=1}^{r-1}(\beta_t+1)\ls m<p_r.$$
This leads a contradiction.

By the above, $\beta_j=\al_j$ for all $j=1,\ldots,r-1$, and thus
$d=n$. We are done. \qed

\medskip

\noindent {\it Proof of Theorem 1.4}.
(i) If $p>2$ is a prime, then $2^{p-1}p=2^{\lfloor(p-1)/(2-1)\rfloor}p$
is a primitive covering number by Theorem 1.3 in the case $r=2$.

By Lemma 2.1, any prime power cannot be a primitive covering number.

Now suppose that $n=p_1^{\al_1} p_2^{\al_2}$ is a primitive covering number,
where $p_1<p_2$ are two distinct primes, and $\al_1,\al_2\in\Z^+$. Then
$$2<\f{\sigma(n)}n<\l(1+\f1{p_1}+\f1{p_1^2}+\cdots\r)\l(1+\f1{p_2}+\f1{p_2^2}+\cdots\r)
=\f{p_1}{p_1-1}\cdot\f{p_2}{p_2-1}.$$
If $p_1>2$, then
$$2<\f{p_1}{p_1-1}\cdot\f{p_2}{p_2-1}\ls\f 3{3-1}\cdot\f 5{5-1}=\f{15}8<2$$
which leads a contradiction. So $p_1=2$. Observe that
$\al_1+1\gs p_2$ by Lemma 2.1. Therefore $n$ is a multiple of $2^{p_2-1}p_2$.
Since both $2^{p_2-1}p_2$ and $n$ are primitive covering numbers, we must have
$n=2^{p_2-1}p_2$.

 (ii) If $p>3$ is a prime, then
 $$2\cdot 3^{\f{p-1}2}p=2^{\f{3-1}{2-1}-1}3^{\f{p-1}{3-1}}p$$
 is a primitive covering number by Theorem 1.3 in the case $r=3$.

 Now assume that $n=p_1^{\al_1} p_2^{\al_2}p_3^{\al_3}$ is a primitive covering number
with $n\eq0\ (\mo\ 3)$,
where $p_1<p_2<p_3$ are distinct primes, and $\al_1,\al_2,\al_3\in\Z^+$.
If $p_1\gs3$, then
$$\align \sum\Sb d\mid n\\d\ \t{is composite}\endSb\f1{\varphi(d)}
<&\sum_{s=1}^3\f1{p_s-1}\l(\f1{p_s}+\f1{p_s^2}+\cdots\r)
\\&+\sum_{1\ls s<t\ls 3}\f1{(p_s-1)(p_t-1)}
\l(1+\f1{p_s}+\f1{p_s^2}+\cdots\r)\l(1+\f1{p_t}+\f1{p_t^2}+\cdots\r)
\\&+\f1{(p_1-1)(p_2-1)(p_3-1)}
\prod_{s=1}^3\l(1+\f1{p_s}+\f1{p_s^2}+\cdots\r)
\\=&\sum_{s=1}^3\f1{(p_s-1)^2}+\sum_{1\ls s<t\ls3}\f{p_sp_t}{(p_s-1)^2(p_t-1)^2}
+\f{p_1p_2p_3}{(p_1-1)^2(p_2-1)^2(p_3-1)^2}
\\\ls&\f1{(3-1)^2}+\f1{(5-1)^2}+\f1{(7-1)^2}+\f{3\cdot 5}{2^24^2}
+\f{3\cdot 7}{2^26^2}+\f{5\cdot 7}{4^26^2}+\f{3\cdot 5\cdot7}{2^24^26^2}=\f{1905}{2304}
\endalign$$
and this contradicts (1.2). So $p_1=2$.
Since $3\mid n$, we have $p_2=3$.
By part (i), $2^2\cdot3$ is a primitive covering number and hence it does not divide $n$.
Therefore $n$ has the form $2\cdot 3^\al p$, where $p>3$ is a prime and $\al\in\Z^+$.

By Lemma 2.1, $(1+1)(\al+1)\gs p$. Thus $\al\gs(p-1)/2$ and hence
$n$ is a multiple of $2\cdot 3^{(p-1)/2}p$. As both $2\cdot 3^{(p-1)/2}p$
and $n$ are primitive covering numbers, we must have $n=2\cdot 3^{(p-1)/2}p$.

(iii) If $p>5$ is a prime, then by Theorem 1.3 both
$$2^35^{\lfloor\f{p-1}4\rfloor}p=2^{\f{5-1}{2-1}-1}5^{\lfloor\f{p-1}{5-1}\rfloor}p
\ \ \ \t{and}\ \ \
2\cdot3\cdot5^{\lfloor\f{p-1}4\rfloor}p=2^{\f{3-1}{2-1}-1}3^{\f{5-1}{3-1}-1}
5^{\lfloor\f{p-1}{5-1}\rfloor}p$$
are primitive covering numbers.

If $p$ is a prime greater than $19$, then
$p>(7-2)(7-3)$, hence both
$$2\cdot3^2\cdot7^{\lfloor\f{p-1}6\rfloor}p=2^{\f{3-1}{2-1}-1}3^{\f{7-1}{3-1}-1}
7^{\lfloor\f{p-1}{7-1}\rfloor}p
\ \ \ \t{and}\ \ \
2^57^{\lfloor\f{p-1}6\rfloor}p=2^{\f{7-1}{2-1}-1}
7^{\lfloor\f{p-1}{7-1}\rfloor}p$$
are primitive covering numbers by Theorem 1.3.
When $p\in\{11,13,17,19\}$, we have
$$p>(7-3)\l(\max\l\{\f{7-1}{3-1},\f{3-1}{2-1}\r\}-1\r)=8$$
and hence $2\cdot3^2\cdot7^{\lfloor\f{p-1}6\rfloor}p$
is still a primitive covering number by the proof of Theorem 1.3.
If $p$ is 11 or 17, then
$$m=(7-1)\l(1-\f{2-1}{7-1}\r)\(\l\lfloor\f{p-1}{7-1}\r\rfloor+1\)<p+1$$
and hence $2^57^{\lfloor\f{p-1}6\rfloor}p$ is a primitive covering number
by the proof of Theorem 1.3.

Combining the above we have shown Theorem 1.4.
\qed

\Ack. The author thanks the referee for helpful comments.
\medskip

\subhead \nofrills{\bf References}\endsubhead

\widestnumber\key{ABCDEF1}

\ref\key E50\by P. Erd\H os\paper On integers of the form $2^k+p$
and some related problems\jour Summa Brasil.
Math.\vol2\yr1950\pages113--123\endref

\ref\key E80\by P. Erd\H os\paper A survey of problems in combinatorial
number theory \jour Ann. Discrete Math.\vol 6\yr 1980\pages 89--115\endref

\ref\key E97\by P. Erd\H os\paper Some of my favorite problems and
results, {\rm in: The mathematics of Paul Erd\H os, I, 47--67,
Algorithms Combin., 13, Springer, Berlin, 1997}\endref

\ref\key F\by M. Filaseta\paper
Coverings of the integers associated with an irreducibility theorem
of A. Schinzel \jour in: Number Theory for the Millennium
(Urbana, IL, 2000), Vol. II, pp. 1-24, A K Peters, Natick, MA,
2002\endref

\ref\key FFKPY\by M. Filaseta, K. Ford, S. Konyagin, C. Pomerance
and G. Yu\paper Sieving by large integers and covering systems
of congruences\jour J. Amer. Math. Soc.
\finalinfo in press. {\tt arXiv:math.NT/0507374}\endref

\ref\key GS\by S. Guo and Z. W. Sun
\paper On odd covering systems with distinct moduli
\jour Adv. in Appl. Math.\vol 35\yr 2005\pages 182--187\endref

\ref\key Gu\by R. K. Guy\book
Unsolved Problems in Number Theory\publ 3rd ed., Springer, New York, 2004\endref

\ref\key H\by J. A. Haight\paper Covering systems of congruences, a negative result
\jour Mathematika\vol 26\yr 1979\pages 53--61\endref

\ref\key PS\by\v S. Porubsk\'y and J. Sch\"onheim \paper Covering
systems of Paul Erd\H os: past, present and future \jour in: Paul
Erd\H os and his Mathematics. I (edited by G. Hal\'asz, L.
Lov\'asz, M. Simonvits, V. T. S\'os), Bolyai Soc. Math. Studies
11, Budapest, 2002, pp. 581--627\endref

\ref\key R\by H. E. Rose\book A Course in Number Theory
\publ Oxford Univ. Press, New York, 1988\endref

\ref \key Si\by R. J. Simpson\paper Regular coverings of the integers by
arithmetic progressions \jour Acta Arith.\vol 45\yr1985\pages
145--152\endref

\ref\key S90\by Z. W. Sun\paper  Finite coverings of groups,
Fund. Math. \vol 134\yr 1990\pages 37--53\endref

\ref\key S96\by Z. W. Sun\paper Covering the integers by arithmetic sequences {\rm II}
\jour Trans. Amer. Math. Soc.\vol348\yr1996\pages4279--4320\endref

\ref\key S00\by Z. W. Sun\paper On integers not of the form $\pm
p^a\pm q^b$ \jour Proc. Amer. Math. Soc.\vol 128\yr 2000\pages
997--1002\endref

\ref\key S01\by Z. W. Sun\paper Algebraic approaches to periodic
arithmetical maps\jour J. Algebra\vol 240\yr
2001\pages723--743\endref

\ref\key S03\by Z. W. Sun\paper Unification of zero-sum problems, subset sums
and covers of $\Z$\jour Electron. Res. Announc. Amer. Math. Soc.
\vol 9\yr 2003\pages 51--60\endref

\ref\key S04\by Z. W. Sun
\paper On the Herzog-Sch\"onheim conjecture for uniform covers of groups \jour
J. Algebra\vol 273\yr 2004\pages153--175\endref

\ref\key S05\by Z. W. Sun\paper On the range of a covering function
\jour J. Number Theory \vol 111\yr 2005\pages 190--196\endref

\ref\key SS\by Z. W. Sun and Z. H. Sun
\paper Some results on covering systems of congruences
\jour J. Southwest-China Teachers Univ.\yr1987\issue1\pages10--15. Zbl. M.
749.11018\endref

\ref\key Z\by \v S. Zn\'am\paper On Mycielski's problem on systems of
arithmetical progressions
\jour Colloq. Math.\vol 15\yr1966\pages 201--204\endref

\enddocument